\documentclass[12pt,a4paper]{article}
\usepackage{amsmath,amssymb,amsfonts}
\def\Bbb{\mathbb}

\title{\bf Dedekind sums in the $p$-adic number field}

\author{Kurt Girstmair}
\date{}

\makeatletter
\let\@@maketitle=\maketitle
\def\maketitle{\def\thispagestyle##1{\relax}\@@maketitle}
\makeatother
\newtheorem{theorem}{Theorem}

\newtheorem{lemma}{Lemma}

\def\BE{\begin{equation}}
\def\EE{\end{equation}}
\def\BD{\begin{displaymath}}
\def\ED{\end{displaymath}}
\def\BA{\begin{array}}
\def\EA{\end{array}}
\def\BEA{\begin{eqnarray*}}
\def\EEA{\end{eqnarray*}}
\def\BI{\bibitem}

\def\N{\Bbb N}
\def\Z{\Bbb Z}
\def\Q{\Bbb Q}
\def\R{\Bbb R}

\def\phi{\varphi}

\def\MB{\mbox}
\def\LD{\ldots}

\def\DIV{\,|\,}
\def\NDIV{\, \nmid \,}

\def\BQ{``}
\def\EQ{'' }

\def\MN{\medskip\noindent}
\def\STOP{\hfill$\Box$}

\def\DED{Dedekind }

\def\QP{\Q_p}
\def\ZP{\Z_p}
\def\ZPT{\Z_p^{\times}}

\begin{document}
\maketitle

\begin{abstract}

\noindent
In a recent note W. Kohnen asks whether the values of Dedekind sums are dense in the field of $p$-adic numbers.
The present paper answers this question. Dedekind sums do not approximate units of $\Z_2$ or $\Z_3$, so they are not dense in $\Q_2$ or $\Q_3$.
But they are dense in $\Q_p$ if $p\ge 5$.
\end{abstract}

\section*{1. Introduction and results}

Let $n$ be a natural number, $m$ an integer, $(m,n)=1$. The classical \DED sum $s(m,n)$ is defined by
\BD
   s(m,n)=\sum_{k=1}^n ((k/n))((mk/n))
\ED
where $((\LD))$ is the \BQ sawtooth function\EQ defined by
\BD
  ((t))=\begin{cases}
                 t-\lfloor t\rfloor-1/2, & \MB{ if } t\in\R\smallsetminus \Z; \\
                 0, & \MB{ if } t\in \Z
       \end{cases}
\ED
(see, for instance, \cite[p. 1]{RaGr}). In the present setting it is more
convenient to work with
\BD
 S(m,n)=12s(m,n)
\ED instead.
Since $S(m+n,n)=S(m,n)$, we obtain all \DED sums if we restrict $m$ to the range $0\le m< n$.

For a prime $p$, let
$\QP$ denote the field of $p$-adic numbers, $\ZP$ the ring of $p$-adic integers, and $\ZPT$ the unit group of $\ZP$.

It is well-known that the values of \DED sums are dense in the field of real numbers (see \cite{Hi,Gi2,Ko}). In the paper \cite{Ko}, W. Kohnen asks
whether the values of \DED sums are dense in $\QP$.
We answer this question as follows.

\begin{theorem} 
\label{t0}

The set of \DED sums is dense neither in $\Q_2$ nor in $\Q_3$, but it is dense in each $\QP$, $p\ge 5$ .

\end{theorem} 

Theorem \ref{t0} is an immediate consequence of Theorems \ref{t1}--\ref{t3} below.
We shall use the following
terminology: Let $a\in \QP$ and let $M$ be a subset of $\QP$. We say that $M$ {\em approximates } $a$ (in $\QP$), if for every integer $k\ge 1$ there is an
element $x\in M$ and an element $b\in\ZP$ such that
\BD
  a=x+p^kb.
\ED
In the case $M=\{S(m,n): n\in\N, 0\le m<n-1, (m,n)=1\}$, we simply say that {\em \DED sums approximate} $a$ (in $\QP$).

\begin{theorem} 
\label{t1}

\DED sums do not approximate any $u\in \ZPT$ for $p=2,3$.

\end{theorem} 

So the answer to Kohnen's question is negative for $p=2,3$. The following result says which $p$-adic integers are approximated by \DED sums.

\begin{theorem} 
\label{t2}

In the cases $p=2,3$, \DED sums approximate each $a\in p\ZP$.\\
If $p\ge 5$, \DED sums approximate each $a\in \ZP$.

\end{theorem} 

Finally, we deal with the approximation of numbers in $\QP\smallsetminus\ZP$.

\begin{theorem} 
\label{t3}
Let $q>1$ be a power of $p$. Then \DED sums approximate each $a\in \frac 1q\ZPT$.

\end{theorem} 

Theorems \ref{t1}--\ref{t3} obviously imply Theorem \ref{t0}.

\section*{2. Proofs}
Theorem \ref{t1} is a consequence of the following lemma. Observe that $nS(m,n)$ is an integer (see \cite[p. 27, Th. 2]{RaGr}).

\begin{lemma} 
\label{l1} Let $n$ be a natural number, $m$ an integer, $(m,n)=1$.\\
{\rm (a)} If $n$ is odd, then $nS(m,n)\equiv 0 \mod 2$.\\
{\rm (b)} If $n\equiv 2 \mod 4$, then $nS(m,n)\equiv 0 \mod 4$.\\
{\rm (c)} If $n\equiv 0 \mod 4$, then $nS(m,n)\equiv 2 \mod 4$.\\
{\rm (d)} If $n\not\equiv 0 \mod 3$, then $nS(m,n)\equiv 0\mod 3$.\\
{\rm (e)} If $n\equiv 0 \mod 3$, then $nS(m,n)\not\equiv 0\mod 3$.
\end{lemma} 

{\em Proof.} Assertions (a) and (d) are immediate consequences of  the aforementioned Theorem 2 in \cite{RaGr} and have been used by various authors
(see, e.g., \cite[formula (69)]{Sa}). Assertion (e) is a weaker form of formula (70) in the said paper \cite{Sa}. Probably assertions (b) and (c)
are also known, but we do not know an appropriate reference. Hence we give a short proof.

Since $n>1$, we may assume that $m$ is in the range $1\le m\le n-1$.
The reciprocity law for \DED sums (see \cite[p. 5]{RaGr})
says
\BD
  mnS(m,n)+ mnS(n,m)=m^2+n^2+1-3mn.
\ED
It is easy to check that the right-hand side is $\equiv 0 \mod 4$ if $n\equiv 2\mod 4$, and $\equiv 2\mod 4$ if $n\equiv 0\mod 4$.
Moreover, $mS(n,m)\equiv 0\mod 2$, by (a). Accordingly, $mnS(n,m)\equiv 0$ mod $4$ if $n$ is even. This shows
\BD
  mnS(m,n)\equiv \left\{\begin{array}{lr}
     0 \mod 4, & \MB{ if }n\equiv 2 \mod 4;\\
     2 \mod 4, & \MB{ if } n\equiv 0 \mod 4.
  \end{array}
  \right.
\ED
Since $m$ is odd, the assertions (b), (c) follow.
\STOP

\MN
{\em Proof of Theorem \ref{t1}.}
First let $u\in \Z_2^{\times}$ and suppose that there is a natural number $n$ and an $m$ with $(m,n)=1$ such that
\BE
\label{2.0}
  u=S(m,n)+4b
\EE
with $b\in \Z_2$.
Then
\BE
\label{2.1}
  nu=nS(m,n)+4nb.
\EE
By assertion (a) of Lemma \ref{l1}, the right-hand side of (\ref{2.1}) is divisible by 2 if $n$ is odd, whereas 2 does not divide $nu$ in this case.
If $n\equiv 2 \mod 4$, assertion (b) says that the right-hand side of (\ref{2.1}) is $\equiv 0\mod 4$, but $nu$ is not divisible by 4.
If $n\equiv 0 \mod 4$, $nu$ and $4nb$ are divisible by 4, whereas $nS(m,n)\equiv 2 \mod 4$, by (c). Hence (\ref{2.1}) and (\ref{2.0}) are impossible.
This settles the case $p=2$ of Theorem \ref{t1}.

The argument in the case $p=3$ is similar. Suppose that
\BE
\label{2.2}
 u=S(m,n)+3b
\EE
holds with $u\in\Z_3^{\times}$ and $b\in\Z_3$.
Then
\BE
\label{2.3}
nu=nS(m,n)+3nb.
\EE
If 3 does not divide $n$, the right-hand side of (\ref{2.3}) is divisible by 3, by assertion (d) of Lemma \ref{l1}.
However, 3 does not divide $nu$. If 3 divides n, then $nu$ and $3nb$ are divisible by 3, whereas
$nS(m,n)\not\equiv 0$ mod 3, by (e). Therefore, (\ref{2.3}) and (\ref{2.2}) are impossible. Hence the theorem holds for $p=3$.
\STOP

\begin{lemma} 
\label{l2}
Let $p$ be a prime. \\
{\rm (a)} Let $q>1$ be a power of $p$, $r\in \Z$ such that $p\NDIV r$, and $s,t\in\Z$. Then the set
\BD
  \{l\in\N:\: l\equiv s\mod q,\: l\equiv t\mod r\}
\ED
approximates each $j\in\Z$, $j\equiv s\mod q$ (in $\QP$). \\
{\rm (b)} If $r$ is as in {\rm (a)} and $t\in\Z$, the set $\{l\in\N; l\equiv t\mod r\}$ approximates each element of $\ZP$.\\
{\rm (c)} If $M\subseteq \ZP$ approximates each element of $\ZP$, so does $uM=\{ux: x\in M\}$, where $u$ is an arbitrary element of $\ZPT$.\\
{\rm (d)} If $M\subseteq \ZP$ approximates each element of $\ZPT$, so does $uM=\{ux: x\in M\}$, where $u$ is an arbitrary element of $\ZPT$.
\end{lemma} 

{\em Proof.} In order to prove (a), let $j\in\Z$, $j\equiv s\mod q$. Let $k\in\N$ be such that $q\DIV p^k$. By the Chinese remainder theorem, there is an $l\in\N$,
$l\equiv j\mod p^k$, $l \equiv t\mod r$. Then $l\equiv s\mod q$ and $j=l+p^kw$ for an integer $w$. Hence (a) follows.
As to (b), let $j, t\in\Z$ be given. By (a), the set $\{l\in\N:l\equiv j\mod p, l\equiv t \mod r\}$ approximates $j$ (in $\QP$). Hence
$\{l\in\N:l\equiv t\mod r\}$ approximates each $j\in\Z$. Since $\Z$ is dense in $\ZP$, assertion (b) follows.
Further, if $M$ approximates each $a\in\ZP$, then $uM$ approximates each $ua$, $a\in \ZP$. Because $u\in\ZPT$, $ua$ takes all values in $\ZP$.
Assertion (d) follows in the same way.
\STOP

\begin{lemma} 
\label{l3}
Let $q>1$ be a power of the prime $p$ and $m\in\N$, $p\NDIV m$.
Put $n=q(m^2+1)$.
Then
\BE
\label{2.5}
  S(m,n)=\frac{q^2-1}q m +S(m,q).
\EE
\end{lemma} 

{\em Proof.} We apply the reciprocity law for \DED sums twice. First,
\BE
\label{2.4}
 S(m,n)=-S(n,m)+\frac nm+\frac mn+\frac 1{nm}-3.
\EE
Since $n\equiv q\mod m$, we have
\BE
\label{2.6}
S(n,m)=S(q,m)=-S(m,q)+\frac qm+\frac mq+\frac 1{qm}-3.
\EE
Inserting the right-hand side of (\ref{2.6}) in (\ref{2.4}) and replacing $n$ by $q(m^2+1)$ yields the assertion.
\rule{5mm}{0mm}\STOP

\MN
Of course, the values of $S(m,q)$ in Lemma \ref{l3} are known for $q=2,3,5$. We note the final form of Lemma \ref{l3} for these cases, which will be needed below.

\begin{lemma} 
\label{l4}
Let $m$ be a natural number. \\
{\rm (a)} If $m$ is odd and $n=2(m^2+1)$, then
\BE
\label{2.8}
  S(m,n)= \frac{3m}2.
\EE
{\rm (b)} If $3\NDIV m$ and $n=3(m^2+1)$, then
\BE
\label{2.10}
   S(m,n)=\begin{cases}
           \frac{2}3(4m+1), & \MB{ if } m\equiv 1\mod 3;\\ \rule{0mm}{5mm}
           \frac{2}3(4m-1), & \MB{ if } m\equiv 2\mod 3.
           \end{cases}
\EE
{\rm (c)} If $5\NDIV m$ and $n=5(m^2+1)$, then
\BE
\label{2.12}
   S(m,n)=\begin{cases}
           \frac{12}5(2m+1), & \MB{ if } m\equiv 1\mod 5;\\ \rule{0mm}{5mm}
           \frac{12}5(2m-1), & \MB{ if } m\equiv 4\mod 5;\\ \rule{0mm}{5mm}
           \frac{24}5m, & \MB{ if } m\equiv 2,3 \mod 5.
           \end{cases}
\EE
\end{lemma} 

\MN
{\em Proof of Theorem \ref{t2}}.
For the time being, let $M$ denote the set of odd natural numbers.
First we deal with the case $p\ge 5$. For $m\in M$ and $n=2(m^2+1)$, (\ref{2.8}) gives
$S(m,n)=3m/2$. Hence \DED sums take all values in $3M/2$.
By Lemma \ref{l2}, (b), the set $M$ approximates each $a\in\ZP$.
Since $2$ and $3$ are in $\ZPT$, Lemma \ref{l2}, (c) says that $3M/2$ approximates each $a\in\ZP$; and so do \DED sums.

If $p=3$, we observe that $M$ approximates each $a \in\Z_3$. Since $2\in\Z_3^{\times}$, $M/2$ also approximates
each $a\in\Z_3$, hence $3M/2$ approximates each $a\in 3\Z_3$. Again, this also holds for the set of \DED sums.

The case $p=2$ is more complicated. First let $l$ be an odd natural number with additional properties specified in the following. Put
\BD
 m=\begin{cases}
           (l-1)/4, & \MB{ if } l\equiv 1\mod 4, \enspace l\equiv 2\mod 3;\\ \rule{0mm}{5mm}
           (l+1)/4, & \MB{ if } l\equiv 3\mod 4, \enspace l\equiv 1\mod 3
           \end{cases}
\ED
and put $n=3(m^2+1)$.
Then (\ref{2.10}) gives $S(m,n)=2l/3$
in both cases.
By Lemma \ref{l2}, (a), the set $\{l\in\N:l\equiv 1 \mod 4, l\equiv 2\mod 3\}$ approximates each $a\in \Z$, $a\equiv 1\mod 4$ (in $\Q_2$).
Further, $\{l\in\N:l\equiv 3 \mod 4, l\equiv 1\mod 3\}$ approximates each $a\in\Z$, $a\equiv 3\mod 4$. Since the set $M$ of odd natural numbers is dense in $\Z_2^{\times}$,
the union $N=\{l\in\N:l\equiv 1 \mod 4, l\equiv 2\mod 3\}\cup\{l\in\N:l\equiv 3 \mod 4, l\equiv 1\mod 3\}$ approximates
each element of $\Z_2^{\times}$. By Lemma \ref{l2}, (d),
$N/3$ approximates each $a$ in $\Z_2^{\times}$. Accordingly, $2N/3$ approximates each $a\in 2\Z_2^{\times}$,
and so do \DED sums.

Next let $l$ be an odd natural number, $l\equiv 3\mod 5$. Then $m=(l-1)/2$ is a natural number $\equiv 1\mod 5$. We put $n=5(m^2+1)$. Now (\ref{2.12}) gives $S(m,n)=12l/5$.
Hence we know that \DED sums take all values $12l/5$, $l\in\N$, $l$ odd,  $l\equiv 3\mod 5$. By Lemma \ref{l2}, (a), these numbers $l$ approximate each odd integer (in $\Q_2$), and, thus,
each element of $\Z_2^{\times}$. Since $3$ and $5$ are in $\Z_2^{\times}$, the numbers $12l/5$ approximate each element of $4\Z_2^{\times}$, and so do \DED sums.

In the final step we use that, by (\ref{2.12}), \DED sums take all values $24m/5$ for $m\in\N$, $m\equiv 2\mod 5$, and $n=5(m^2+1)$ (observe that $m$ need not be odd). By our above arguments,
\DED sums approximate each $a\in 8\Z_2$. This concludes the proof.
\STOP

\begin{lemma} 
\label{l5}
Let $q>1$ be a power of the prime $p$, $r\in \{1,\LD,q-1\}$, $p\NDIV r$, and let $r^*\in\{1,\LD,q-1\}$ be defined by $rr^*\equiv 1\mod q$.
Then \DED sums take all values $l/q$, where $l$ is a natural number and
\BE
\label{2.15}
  l\equiv r^*\mod q,\enspace l\equiv qS(r,q)\mod q^2-1.
\EE
\end{lemma} 

{\em Proof.} First we note
\BE
\label{2.16}
 qS(r,q)\equiv  r+r^* \mod q.
\EE
This is well-known, but for the sake of convenience we give a short proof.
The reciprocity law yields
\BD
rqS(r,q)=-rqS(q,r)+r^2+q^2+1-3rq.
\ED
Now $rS(q,r)$ is an integer, and so $rqS(q,r)\equiv 0\mod q$. Accordingly,
\BD
 rqS(r,q)\equiv r^2+1\mod q,
\ED
and from $rr^*\equiv 1\mod q$ we obtain (\ref{2.16}).

Let $m$ be a natural number, $m\equiv r\mod q$, and put $n=q(m^2+1)$. Then (\ref{2.5}) says
\BE
\label{2.18}
 S(m,n)=\frac{q^2-1}q m+S(r,q).
\EE
Accordingly, $qS(m,n)$ is an integer, and $qS(m,n)\equiv (q^2-1)m+qS(r,q)\equiv -m+r+r^*\mod q$, by (\ref{2.16}). Since $m\equiv r\mod q$,
we obtain $qS(m,n)\equiv r^*\mod q$. On the other hand, (\ref{2.18}) yields $qS(m,n)\equiv qS(r,q)\mod q^2-1$. Hence $S(m,n)$ takes the form
$S(m,n)=l/q$, where the integer $l$ has the properties of (\ref{2.15}).

Conversely, we show that for each natural number $l$ with the properties of (\ref{2.15}) we obtain $l/q$ as the value of a \DED sum. To this
end we put
\BE
\label{2.22}
  m=\frac{l-qS(r,q)}{q^2-1}.
\EE
Since $l\equiv qS(r,q)\mod q^2-1$, $m$ is an integer. Further,
\BD
(q^2-1)m\equiv -m \equiv l-qS(r,q)\equiv r^*-(r+r^*)\mod q,
\ED
by (\ref{2.16}), and so $m\equiv r\mod q$. If $S(r,q)\le 0$, then $m$ is positive. If $S(r,q)>0$, we observe
\BD
 S(r,q)\le S(1,q)=q-3+\frac 2q\le q-2
\ED
(see \cite[Satz 2]{Ra}). Therefore, $0<qS(r,q)<q^2-1$. Since the natural number $l$ is $\equiv qS(r,q)\mod q^2-1$, $l$ must be $\ge qS(r,q)$. However, $l=qS(r,q)$ is impossible,
since $l-qS(r,q)\equiv -r\mod q$, and $(r,q)=1$. Accordingly, $m$ is a natural number, $m\equiv r\mod q$. If we put $n=q(m^2+1)$, we obtain
(\ref{2.18}), which, by (\ref{2.22}), is reduced to $S(m,n)=l/q$.
\STOP

\MN
{\em Proof of Theorem \ref{t3}}. Let $r\in \Z$, $p\NDIV r$. By Lemma \ref{l5}, \DED sums take all values $l/q$, $l\in\N$, $l\equiv r^*\mod q$,
$l\equiv qS(r,q)\mod q^2-1$. The set of these numbers $l$, however, approximates each number $j/q$, $j\in\Z$, $j\equiv r^*\mod q$ (in $\QP$), by Lemma \ref{l2}, (a).
If we vary $r\in\{1,\LD,q\}$, $p\NDIV r$, we see that \DED sums
approximate all numbers $j/q$, $j\in \Z$, $p\NDIV j$ (in $\QP$). But the set of these numbers $j/q$ is dense in $\frac 1q\ZPT$.
\STOP


\vspace{0.5cm}
\noindent
Kurt Girstmair            \\
Institut f\"ur Mathematik \\
Universit\"at Innsbruck   \\
Technikerstr. 13/7        \\
A-6020 Innsbruck, Austria \\
Kurt.Girstmair@uibk.ac.at

\end{document}